\documentclass[11pt]{article}
\usepackage{graphicx}
\usepackage{amsmath}
\usepackage{epsfig}  
\usepackage{verbatim}  
\usepackage{tabularx}
\usepackage{epic, eepic}
\usepackage{xcolor}
\usepackage{lineno}
\usepackage{amssymb,latexsym,amscd,amsmath,amsfonts,enumerate,supertabular}
\usepackage{subfigure}
\usepackage{longtable}
\newtheorem{Theorem}{Theorem}[section] 

\newtheorem{lemma}[Theorem]{Lemma} 
\newtheorem{Proof}[Theorem]{Proof}
\newtheorem{Proposition}[Theorem]{Proposition}

\begin{document}

\title{{\bf A simple numerical method of second and third orders convergence for solving a fully third order nonlinear boundary value problem }}

\author{ Dang Quang A$^{\text a}$,  Dang Quang Long$^{\text b}$\\
$^{\text a}$ {\it\small Center for Informatics and Computing, VAST}\\
{\it\small 18 Hoang Quoc Viet, Cau Giay, Hanoi, Vietnam}\\
{\small Email: dangquanga@cic.vast.vn}\\
$^{\text b}$ {\it\small Institute of Information Technology, VAST,}\\
{\it\small 18 Hoang Quoc Viet, Cau Giay, Hanoi, Vietnam}\\
{\small Email: dqlong88@gmail.com}}
\date{ }
%
\date{}          
\maketitle
\begin{abstract}
In this paper we consider a fully third order nonlinear boundary value problem which is of great interest of many researchers. First we establish the existence, uniqueness  of solution. Next, we propose  simple iterative methods on both continuous and discrete levels. We prove that the discrete methods are of second order and third accuracy due to the use of appropriate formulas for numerical integration and obtain estimate for total error.  Some examples demonstrate the validity of the obtained theoretical results and the efficiency of the iterative method.
\end{abstract}
{\small
}
\section{Introduction}
Boundary value problems (BVPs) for third order nonlinear differential equations  appear in many applied fields, such as flexibility mechanics, chemical engineering, heat conduction and so on. A lot of works are devoted to the qualitative aspects of the problems (see e.g. \cite{Bai,Cabada2,Gross,Guo1,Rez,Sun,Zhai}). There are also many methods concerning the solution of third order BVPs including analytical methods \cite{Abus,Lv-Gao,Pue} and numerical methods by using interpolation polynomials \cite{Al-Sa}, quartic splines \cite{Gao}, \cite{Pandey2}, quintic splines\cite{Khan}, Non-polynomial splines  \cite{Isla1}, \cite{Isla2}, \cite{Sriv}, and wavelet \cite{Faza}. The majority of the mentioned above numerical methods are devoted to linear  equations or special nonlinear third order differential equations.\\
In this paper we consider the following BVP
\begin{equation}\label{eq1}
\begin{split}
u^{(3)}(t) &=f(t, u(t), u'(t), u''(t)), \quad 0 < t < 1, \\
u(0)&=c_1, u'(0)=c_2, u'(1)=c_3.
\end{split}
\end{equation}
Some authors studied the existence and positivity of solution for this problem, for example, by using the lower and upper solutions method and fixed point theorem on cones, in \cite{Yao-Feng} Yao and Feng established the existence of solution and positive solution for the case $f=f(t,u(t))$, in \cite{Feng-Liu} Feng and Liu obtained existence results by the use of the lower and upper solutions method and a new maximum principle for the case $f=f(t,u(t), u'(t))$. It should be emphasized that the results of these two works are pure existence but not methods for finding solutions. Many researchers are interested in numerical solution of the problem \eqref{eq1} without attention to qualitative aspects of it or refer to the book \cite{Agarwal}. \par

Below we mention some works devoted to solution methods for the problem \eqref{eq1}. Namely, 
Al Said et al. \cite{Al-Said} have solved a  third order two point BVP using cubic splines. Noor 
et al. \cite{Noor} generated second order method based on quartic splines. Other authors \cite{Cala,Khan} generated finite difference using fourth degree B-spline and quintic polynomial spline for this problem subject to other boundary conditions. El-Danaf \cite{Danaf} constructed a new spline method based on quartic nonpolynomial spline functions that has a polynomial part and a trigonometric part to develop numerical methods for a linear differential equation with the boundary conditions as in \eqref{eq1}. Recently, in 2016 Pandey \cite{Pandey1} solved the problem for the case $f=f(t,u)$ by the use of quartic polynomial splines. The convergence  of the method at least $O(h^2)$ for the linear case $f=f(t)$ was proved.  In the next year this author in \cite{Pandey2} proposed two difference schemes for the general case 
$f=f(t, u(t), u'(t), u''(t))$ and also established the second order accuracy for the linear case. In the beginning of 2019 Chaurasia et al. \cite{Chau} use exponential amalgamation of cubic spline functions to form a novel numerical method of second-order accuracy. 
   It should be emphasized that all mentioned above authors only draw attention to the construction of discrete analog of the problem \eqref{eq1} and estimate the error of the obtained  solution assuming that the nonlinear system of algebraic equations can be solved by  known iterative methods. Thus, they did not take into account the errors arisen in the last iterative methods.\par
Motivated by these facts, in this paper we propose a completely different method, specifically, 
an iterative method on both continuous and discrete levels for the problem \eqref{eq1}. We give an analysis of total error of the solution actually obtained. This error includes the error of the iterative method on continuous level and the error arisen in numerical realization of this iterative method. The obtained total error estimate suggests to choose suitable grid size for discretization if desiring to get approximate solution with a given accuracy. In order to justify the total error estimate, first we establish some results on existence, uniqueness of solution. These results are obtained by the method developed in \cite{Dang1}-\cite{Dang8}.  Some examples demonstrate the validity of the obtained theoretical results and the efficiency of the iterative method.

\section{Existence results}
For simplicity of presentation we consider the problem \eqref{eq1} with homogeneous boundary conditions, i.e., the problem
\begin{equation}\label{eq2}
\begin{split}
u^{(3)}(t) &=f(t, u(t), u'(t), u''(t)), \quad 0 < t < 1, \\
u(0)&=0, u'(0)=0, u'(1)=0.
\end{split}
\end{equation}

To investigate this problem  we associate it with an operator equation as follows.\par 
For functions $\varphi (x) \in C[0, 1]$ consider the nonlinear operator $A$ defined by
\begin{equation}\label{defA}
(A\varphi )(t)=f(t,u(t), u'(t), u''(t)),
\end{equation}
where $u(t)$ is the solution of the problem
\begin{equation}\label{eq3}
\begin{split}
u'''(t)&=\varphi (t), \quad 0<t<1\\
u(0)&=0, u'(0)=0, u'(1)=0.
\end{split}
\end{equation}

\begin{Proposition}\label{prop1}   If the function $\varphi (x)$ is a fixed point of the operator $A$, i.e., $\varphi (t)$ is a solution of the operator equation
$\varphi = A\varphi $ , 
then the function $u(t)$ determined from the BVP \eqref{eq3} solves the problem \eqref{eq2}. Conversely, if $u(t)$ is a solution of the BVP \eqref{eq2} then the function
$\varphi(t)=f(t,u(t), u'(t), u''(t))$
is a fixed point of the operator $A$ defined above by \eqref{defA}, \eqref{eq3}.
\end{Proposition}
Thus,  the problem \eqref{eq2} is reduced to the fixed point problem for $A$.\par 
Now, we study the properties of $A$. For this purpose, notice that the problem \eqref{eq2}  has a unique solution representable in the form
\begin{equation}\label{eq2.8}
u(t)=\int_0^1 G_0(t,s)\varphi(s)ds, \quad 0<t<1,
\end{equation}
where $G_0(t,s)$ is the Green function of the problem \eqref{eq3}
\begin{equation*}
\begin{aligned}
G_0(t,s)=\left\{\begin{array}{ll}
\dfrac{s}{2}(t^2-2t+s), \quad 0\le s \le t \le 1,\\
\, \, \dfrac{t^2}{2}(s-1), \quad 0\le t \le s \le 1.\\
\end{array}\right.
\end{aligned}
\end{equation*}
Differentiating both sides of \eqref{eq2.8} gives
\begin{align}
u'(t) & =\int_0^1 G_1(t,s)\varphi(s)ds, \label{eq2.8a}\\
u''(t) & =\int_0^1 G_2(t,s)\varphi(s)ds, \label{eq2.8b}
\end{align}
where
\begin{equation*}
G_1(t,s)=\left\{\begin{array}{ll}
s(t-1), \quad 0\le s \le t \le 1,\\
t(s-1), \quad 0\le t \le s \le 1,\\
\end{array}\right.
\end{equation*}
\begin{equation}\label{eqG2}
G_2(t,s)=\left\{\begin{array}{ll}
s,& \quad 0\le s \le t \le 1,\\
s-1,& \quad 0\le t \le s \le 1.\\
\end{array}\right.
\end{equation}

It is easily seen that
$G_0(t,s) \leq 0, \; G_1(t,s) \leq 0$
in $Q=[0,1]^2$ and
\begin{equation}\label{valGreen}
\begin{aligned}
&M_0= \max _{0\le t\le 1} \int _0 ^1 |G(t,s)| \ ds =\dfrac{1}{12},\quad M_1= \max _{0\le t\le 1} \int _0 ^1 |G_1(t,s)| \ ds =\dfrac{1}{8},\\
&M_2= \max _{0\le t\le 1} \int _0 ^1  |G_2(t,s)|  \ ds =\dfrac{1}{2}.\\
\end{aligned}
\end{equation}

Next, for each fixed real number $M>0$ introduce the domain
\begin{equation*}
\mathcal{D}_M=\{ (t,x,y,z)| \ 0\leq t\leq 1, \,\,|x| \leq M_0M, \,\, |y| \leq M_1M, \,\, 
 |z| \leq M_2M \},
\end{equation*}
and as usual, by $B[O,M]$ we denote the closed ball of radius $M$ centered at $0$ in the space of continuous in $[0, 1]$ functions, namely,
$
B[O,M]=\{ \varphi \in C[0,1]| \ \| \varphi \| \leq M \},
$
where
$\| \varphi \|=  \max_{0 \leq t \leq 1} |\varphi (t)|. $

By the analogous techniques as in \cite{Dang1}-\cite{Dang8} we have proved the following results.
\begin{Theorem}[Existence of solutions]\label{theorem1}
Suppose that there exists a number $M>0$ such that the function $f(t,x,y,z)$ is continuous and bounded by $M$ in the domain $\mathcal{D}_M$, i.e., 
\begin{equation*}
|f(t,x,y,z)| \leq M
\end{equation*}
for any $(t,x,y,z) \in \mathcal{D}_M .$\par
Then, the problem \eqref{eq1} has a solution $u(t)$ satisfying 
\begin{equation*}
|u(t)| \leq M_0M, \; |u'(t)| \leq M_1M,\; |u''(t)| \leq M_2M \text{ for any } 0 \le t \le 1.
\end{equation*}
\end{Theorem}
\begin{Theorem}[Existence and uniqueness of solution]\label{theorem3}
Assume that there exist numbers
  $M,L_0, L_1$,  $ L_2 \geq 0$ such that
\begin{equation*}
|f(t,x,y,z)| \leq M,
\end{equation*}
\begin{multline}
|f(t,x_2,y_2,z_2)-f(t,x_1,y_1,z_1)| \leq  L_0|x_2-x_1|+ L_1|y_2-y_1|+L_2|z_2-z_1|       
\end{multline}
for any $(t,x,y,z), (t,x_i,y_i,z_i) \in \mathcal{D}_M \ (i=1,2)$ and
\begin{equation*}
q:=L_0M_0+ L_1M_1+L_2M_2<1.
\end{equation*}
Then, the problem \eqref{eq2} has a unique solution $u(t)$ such that $|u(t)| \leq M_0M,$ $|u'(t)| \leq M_1M, \,\, |u''(t)| \leq M_2M $ for any $0 \le t \le 1$.
\end{Theorem}
\noindent {\bf Remark.} The problem \eqref{eq1} for $u(t)$ with non-homogeneous boundary conditions can be reduced to the problem with homogeneous for function $v(t)$ if setting $u(t)=v(t)+P_2(t)$, where $P_2(t)$ is the second degree polynomial satisfying the boundary conditions $P_2(0)=c_1, P'_2(0)=c_2, P_2(1)=c_3$.

\section{Iterative method on continuous level}\label{IterMeth}
Consider the following iterative method for solving the problem \eqref{eq2}:
\begin{enumerate}
\item Given 
\begin{equation}\label{iter1c}
\varphi_0(t)=f(t,0,0,0).
\end{equation}
\item Knowing $\varphi_k(t)$  $(k=0,1,...)$ compute
\begin{equation}\label{iter2c}
\begin{split}
u_k(t) &= \int_0^1 G_0(t,s)\varphi_k(s)ds ,\\
y_k(t) &= \int_0^1 G_1(t,s)\varphi_k(s)ds ,\\
z_k(t) &= \int_0^1 G_2(t,s)\varphi_k(s)ds ,
\end{split}
\end{equation}
\item Update
\begin{equation}\label{iter3c}
\varphi_{k+1}(t) = f(t,u_k(t),y_k(t),z_k(t)).
\end{equation}
\end{enumerate}
Set 
\begin{equation*}
p_k=\dfrac{q^k}{1-q}\| \varphi _1 -\varphi _0\|.
\end{equation*}
\begin{Theorem}[Convergence]\label{theorem5} Under the assumptions of Theorem \ref{theorem3} the above iterative method converges and there hold the estimates
\begin{equation*}
\|u_k-u\| \leq M_0p_k, \quad  \|u'_k-u'\| \leq M_1p_k,\quad
\|u''_k-u''\| \leq M_2p_k,
\end{equation*}
where $u$ is the exact solution of the problem \eqref{eq2} and $M_0, M_1, M_2$ are given by \eqref{valGreen}.
\end{Theorem}
This theorem follows straightforward from the convergence of the successive approximation method for finding fixed point of the operator $A$ and the representations \eqref{eq2.8}-\eqref{eq2.8b} and \eqref{iter2c}.
\section{Discrete iterative method 1}
To numerically realize the above iterative method we construct the corresponding discrete iterative methods. For this purpose cover the interval $[0, 1]$   by the uniform grid $\bar{\omega}_h=\{t_i=ih, \; h=1/N, i=0,1,...,N  \}$ and denote by $\Phi_k(t), U_k(t), Y_k(t),  Z_k(t)$ the grid functions, which are defined on the grid $\bar{\omega}_h$ and approximate the functions $\varphi_k (t), u_k(t), y_k(t),  z_k(t)$ on this grid, respectively.\par
First, consider the following discrete iterative method, named {\bf Method 1}:
\begin{enumerate}
\item Given 
\begin{equation}\label{iter1d}
\Phi_0(t_i)=f(t_i,0,0,0),\ i=0,...,N. 
\end{equation}
\item Knowing $\Phi_k(t_i),\; k=0,1,...; \; i=0,...,N, $  compute approximately the definite integrals \eqref{iter2c} by trapezium formulas
\begin{equation}\label{iter2d}
\begin{split}
U_k(t_i) &= \sum _{j=0}^N h\rho_j G_0(t_i,t_j)\Phi_k(t_j), \\
Y_k(t_i) &= \sum _{j=0}^N h\rho_j G_1(t_i,t_j)\Phi_k(t_j) ,\\
Z_k(t_i) &= \sum _{j=0}^N h\rho_j G_2^*(t_i,t_j)\Phi_k(t_j) ,\;  i=0,...,N,
\end{split}
\end{equation}
\noindent where $\rho_j$ are the weights of the trapezium formula
\begin{equation*}
\rho_j = 
\begin{cases}
1/2,\; j=0,N\\
1, \; j=1,2,...,N-1
\end{cases}
\end{equation*}
and 
\begin{equation}\label{eqg2*}
G_2^*(t,s)  =
\begin{cases}
s, \quad & 0\leq s < t\leq 1, \\
s-1/2, \quad & s=t, \\
s-1, & 0\leq t < s\leq 1.
\end{cases}
\end{equation}
\item Update
\begin{equation}\label{iter3d}
\Phi_{k+1}(t_i) = f(t_i,U_k(t_i),Y_k(t_i),Z_k(t_i)).
\end{equation}
\end{enumerate}
In order to get the error estimates for the numerical approximate solution for $u(t)$ and its derivatives on the grid we need some following auxiliary results.
\begin{Proposition}\label{prop1}
Assume that the function $f(t,x,y,z)$ has all continuous partial derivatives up to second order in the domain $\mathcal{D}_M$. Then for the functions $u_k(t), y_k(t),  z_k(t), k=0,1,...$, constructed by the iterative method \eqref{iter1c}-\eqref{iter3c} there hold
$z_k(t) \in C^3 [0, 1], \;  y_k(t) \in C^4 [0, 1], \;u_k(t) \in C^5 [0, 1].$
\end{Proposition}
\begin{Proof}
We prove the proposition by induction. For $k=0,$ by the assumption on the function $f$ we have $\varphi_0(t) \in C^2[0, 1]$ since $\varphi_0(t)=f(t,0,0,0)$. Taking into account the expression \eqref{eqG2} of the function $G_2(t,s)$ we have
\begin{equation*}
z_0(t)=\int _0^1 G_2(t,s) \varphi_0(s) ds=
\int_0^t s \varphi_0(s) ds -\int_t^1 (s-1) \varphi_0(s) ds.
\end{equation*}
It is easy to see that $z_0'(t)=\varphi_0(t)$. Therefore, $z_0(t) \in C^3[0, 1]$. It implies $y_0(t) \in C^4 [0, 1], \;u_0(t) \in C^5 [0, 1]$.\par 
Now suppose $z_k(t) \in C^3 [0, 1], \; y_k(t) \in C^4 [0, 1], \;u_k(t) \in C^5 [0, 1].$ Then, because
$\varphi_{k+1}(t) = f(t,u_k(t),y_k(t),z_k(t))$
and the function $f$ by the assumption has continuous derivative in all variables up to order 2, it follows that $\varphi_{k+1}(t) \in C^2[0, 1]$. Repeating the same argument as for $\varphi_0(t)$ above we obtain that $z_{k+1}(t) \in C^3 [0, 1], \; y_{k+1}(t) \in C^4 [0, 1], \;u_{k+1}(t) \in C^5 [0, 1].$
Thus, the proposition is proved.
\end{Proof}

\begin{Proposition}\label{prop2}

For any function $\varphi (t) \in C^2[0, 1]$ there hold the estimates
\begin{equation}\label{eq:prop2}
\int_0^1 G_n (t_i,s) \varphi (s) ds = \sum _{j=0}^N h\rho_j G_n(t_i,t_j)\varphi(t_j) +O(h^2), 
\quad (n=0,1)
\end{equation}
\begin{equation}\label{eq:prop2a}
\int_0^1 G_2 (t_i,s) \varphi (s) ds = \sum _{j=0}^N h\rho_j G_2^*(t_i,t_j)\varphi(t_j) +O(h^2). 
\end{equation}
\end{Proposition}
\begin{Proof}
In the case $n=0, 1$,  since the functions $G_n(t_i,s)$ are continuous at $s=t_i$ and are polynomials in $s$ in the intervals $[0, t_i]$ and $[t_i, 1]$  we have
\begin{align*}
& \int_0^1 G_n (t_i,s) \varphi (s) ds =\int_0^{t_i} G_n (t_i,s) \varphi (s) ds + \int_{t_i}^1 G_n (t_i,s) \varphi (s) ds\\
&= h\big( \tfrac{1}{2}G_n (t_i,t_0)\varphi (t_0)+ G_n (t_i,t_1)\varphi (t_1)+...+G_n (t_i,t_{i-1})\varphi (t_{i-1})+ \tfrac{1}{2}G_2 (t_i,t_i)\varphi (t_i)  \big)\\
& + h\big( \tfrac{1}{2}G_n (t_i,t_i)\varphi (t_i)+ G_n (t_i,t_{i+1})\varphi (t_{i+1})+...+G_n (t_i,t_{N-1})\varphi (t_{N-1})\\
&+ \tfrac{1}{2}G_n (t_i,t_N)\varphi (t_N)  \big) +O(h^2)\\
&= \sum _{j=0}^N h\rho_j G_n(t_i,t_j)\varphi(t_j) +O(h^2)\quad (n=0,1).
\end{align*}
Thus, the estimate \eqref{eq:prop2} is established.
The estimate \eqref{eq:prop2a} is  obtained using the following result, which is easily proved.
\begin{lemma}\label{le2}
Let $p(t)$ be a function having continuous derivatives up to second order in the interval $[0,1]$ except for the point $t_i, \ 0<t_i<1$, where it has a jump.
 Denote $\lim _{t \rightarrow {t_i-0}}p(t)=p_i^- $, 
$\lim _{t \rightarrow {t_i+0}}p(t)=p_i^+,$ 
$   p_i= \tfrac{1}{2}(p_i^- +p_i^+) $. 
Then
\begin{equation}
\int_0^1 p(t) dt= \sum _{j=0}^N h\rho_j p(j) +O(h^2),
\end{equation}
where $p_j=p(t_j), j \ne i.$
\end{lemma}
\end{Proof}

\begin{Proposition}\label{prop3}
Under the assumption of Proposition \ref{prop1}  for any $k=0,1,...$ there hold the estimates
\begin{equation}\label{eq:prop3a}
\|\Phi_k -\varphi_k  \|= O(h^2),
\end{equation}
\begin{equation}\label{eq:prop3b}
\begin{split}
\|U_k -u_k  \|&=O(h^2), \; \|Y_k -y_k  \|=O(h^2),  \; \|Z_k -z_k  \|=O(h^2).
\end{split}
\end{equation}
\noindent where $\|.\|_{C(\bar{\omega}_h)}$ is the max-norm of function on the grid $\bar{\omega}_h$.
\end{Proposition}

\begin{Proof}
We prove the proposition by induction. For $k=0$ we have immediately $\|\Phi_0 -\varphi_0  \|= 0$. Next, by the first equation in \eqref{iter2c} and Proposition \ref{prop2} we have
\begin{equation}
u_0(t_i)=\int_0^1 G_0 (t_i,s) \varphi_0 (s) ds = \sum _{j=0}^N h\rho_j G_0(t_i,t_j)\varphi_0(t_j)+O(h^2)
\end{equation}
for any $i=0,...,N$.
On the other hand, in view of  the first equation in \eqref{iter2d}  we have
\begin{equation}
U_0(t_i)= \sum _{j=0}^N h\rho_j G_0(t_i,t_j)\varphi_0(t_j).
\end{equation}
Therefore, $|U_0(t_i)- u_0(t_i)|= O(h^2)$. Consequently, $\|U_0 -u_0  \|=O(h^2) $.\\
Similarly, we have 
\begin{equation}
\|Y_0 -y_0  \|=O(h^2), \;  \|Z_0 -z_0  \|=O(h^2).
\end{equation}
Now suppose that \eqref{eq:prop3a} and \eqref{eq:prop3b} are valid for $k \ge 0$. We shall show that these estimates are valid for $k+1$.\\
Indeed, by the Lipshitz condition of the function $f$ and the estimates  \eqref{eq:prop3b} it is easy to obtain the estimate
\begin{equation}\label{eq:diffPhi}
\|\Phi_{k+1} -\varphi_{k+1}  \|= O(h^2)
\end{equation}
Now from the first equation in \eqref{iter2c} by Proposition \ref{prop2} we have
\begin{equation*}
u_{k+1}(t_i)=\int_0^1 G_0 (t_i,s) \varphi_{k+1} (s) ds = \sum _{j=0}^N h\rho_j G_0(t_i,t_j)\varphi_{k+1}(t_j)+O(h^2)
\end{equation*}
On the other hand by the first formula in \eqref{iter2d} we have
\begin{equation*}
U_{k+1}(t_i) = \sum _{j=0}^N h\rho_j G_0(t_i,t_j)\Phi_{k+1}(t_j).
\end{equation*}
From the above equalities, 
having in mind the estimate \eqref{eq:diffPhi} we obtain the estimate
$$\|U_{k+1} -u_{k+1}  \|=O(h^2).
$$
Similarly, we obtain
$$\|Y_{k+1} -y_{k+1} \|=O(h^2),\;  \|Z_{k+1} -z_{k+1}\|=O(h^2).
$$
Thus, by induction we have proved the proposition.
\end{Proof}
Now combining Proposition \ref{prop3} and Theorem \ref{theorem5} results in the following theorem.
\begin{Theorem}\label{thm4}
For the approximate solution of the problem \eqref{eq2} obtained by the discrete iterative method on the uniform grid with gridsize $h$ there hold the estimates
\begin{equation*}
\begin{split}
\|U_k-u\| &\leq \left(M_0+\frac{1}{r}\right)p_kd +O(h^2), \; \|Y_k-u'\| \leq M_1p_kd +O(h^2), \\
\|Z_k-u''\| &\leq M_2p_kd +O(h^2).
\end{split}
\end{equation*}
\end{Theorem}

\noindent {\bf Remark 1.} We perform the discrete iterative process \eqref{iter1d}-\eqref{iter3d} until $ \|\Phi_{k+1}-\Phi_k \| \le TOL$, where $TOL$ is a given tolerance. From Theorem \ref{thm4} it is seen that the accuracy of the discrete approximate solution depends on both the number $q$ defined in Theorem \ref{theorem3}, which determines the number of iterations of the continuous iterative method and the gridsize $h$.  
The number $q$ presents the nature of the BVP, therefore, it is necessary to choose appropriate $h$ consistent with $q$ because the choice of very small $h$ does not increase the accuracy of the approximate discrete solution. Below, in examples we shall see this fact.\\


\noindent {\bf Remark 2.} As mentioned in the Introduction, in 2016 Pandey \cite{Pandey1} discretized the problem \eqref{eq1} by quartic splines and proved the second order convergence only for the linear case (when $f=f(x)$). Next year, in \cite{Pandey2} he constructed two difference schemes for the problem and also proved the second order convergence for the linear case. The obtained system of difference equations are solve iteratively by the Gauss-Seidel or Newton-Raphson method. The error arising in these iterative methods are not considered together with the error of the discretization. 

\section{Discrete iterative method 2}
Consider another discrete iterative method, named {\bf Method 2}. The steps of this method are the same of Method 1 with essential difference in Step 2 and now the number of grid points is even, $N=2n$. Namely,\\

2'. Knowing $\Phi_k(t_i),\; k=0,1,...; \; i=0,...,N, $  compute approximately the definite integrals \eqref{iter2c} by the modified  Simpson formulas
\begin{equation}\label{iter2dNew}
\begin{split}
U_k(t_i) &= F(G_0 (t_i,.)\Phi_k(.)),\\
Y_k(t_i) &= F(G_1 (t_i,.)\Phi_k(.)),\\
Z_k(t_i) &= F(G_2^* (t_i,.)\Phi_k(.)),
\end{split}
\end{equation}
where
\begin{equation}\label{•}
F(G_l (t_i,.)\Phi_k(.)) = 
\begin{cases}
\sum _{j=0}^N h\rho_j G_l(t_i,t_j)\Phi_k(t_j) \; \text{ if } i \text {  is  even }\\
\sum _{j=0}^N h\rho_j G_l(t_i,t_j)\Phi_k(t_j) +\dfrac{h}{6}\Big ( G_l(t_i, t_{i-1})\Phi_k(t_{i-1}) -2G_l(t_{i}, t_i)\Phi_k(t_i)\\
\quad +G_l(t_{i} , t_{i+1})\Phi_k(t_{i+1})  \Big )  \; 
   \text{ if } i \text {  is  odd } ,               \\
l=0,1; \; i=0,1,2,...,N.
\end{cases}
\end{equation}
\noindent $\rho_j$ are the weights of the Simpson formula
\begin{equation*}
\rho_j = 
\begin{cases}
1/3,\; j=0,N\\
4/3, \; j=1,3,...,N-1\\
2/3, \; j=2,4,...,N-2,
\end{cases}
\end{equation*}
$F(G_2^* (t_i,.)\Phi_k(.))$ is calculated in the same way as  $F(G_l (t_i,.)\Phi_k(.))$ above, where $G_l$ is replaced by $G_2^*$ defined by the formula \eqref{eqg2*}.

\begin{Proposition}\label{prop1New}
Assume that the function $f(t,x,y,z)$ has all continuous partial derivatives up to fourth order in the domain $\mathcal{D}_M$. Then for the functions $u_k(t), y_k(t),  z_k(t), \varphi_{k+1}(t)$,  $k=0,1,...$, constructed by the iterative method \eqref{iter1c}-\eqref{iter3c} there hold
$z_k(t) \in C^5 [0, 1], \;  y_k(t) \in C^6 [0, 1], \;u_k(t) \in C^7 [0, 1], \varphi_{k+1}(t) \in C^4 [0,1].$
\end{Proposition}

\begin{Proposition}\label{prop2New}

For any function $\varphi (t) \in C^4[0, 1]$ there hold the estimates
\begin{equation}\label{eq:Sim1}
\int_0^1 G_l (t_i,s) \varphi (s) ds = F(G_l (t_i,.)\varphi(.)) +O(h^3), 
\quad (l=0,1)
\end{equation}
\begin{equation}\label{eq:Sim2}
\int_0^1 G_2 (t_i,s) \varphi (s) ds = F(G_2^* (t_i,.)\varphi(.)) +O(h^3). 
\end{equation}
\end{Proposition}
\begin{Proof}
Recall that the interval $[0,1]$ is divided into $N=2n $ by the points $t_i =ih, h=1/N$. In each subinterval $[0, t_i]$ and $[t_i, 1$ the functions $G_l(t_i,s)$ are continuous as polynomials. Therefore, if $i$ is even number, $i=2m$ then we represent
\begin{equation*}
\int_0^1 G_l (t_i,s) \varphi (s) ds=\int _0^{t_{2m}}  \; + \int_{t_{2m}}^1 .
\end{equation*}
Applying the Simpson formula to the integrals in the right-hand side we obtain
\begin{equation*}
\int_0^1 G_l (t_i,s) \varphi (s) ds= F(G_l (t_i,.)\varphi(.)) +O(h^4)
\end{equation*}
because by assumption $\varphi (t) \in C^4[0, 1]$.\\
Now consider the case when $i$ is odd number, $i=2m+1$. In this case we represent
\begin{equation}\label{eqS3}
I= \int_0^1 G_l (t_i,s) \varphi (s) ds=\int _0^{t_{2m}}  \; + \int _{t_{2m}}^{t_{2m+1}} + \int _{t_{2m+1}}^{t_{2m+2}}+
\int_{t_{2m+2}}^1 .
\end{equation}
For simplicity we denote
$$f_j=G_l (t_i,s_j) \varphi (s_j)
$$
Applying the Simpson formula to the first and the fourth integrals in the right-hand side \eqref{eqS3} and the trapezium formula to the second and the third integrals there, we obtain
\begin{align*}
I &=\dfrac{h}{3}[ f_0+f_{2m}+4(f_1+f_3+...+f_{2m-1})+2(f_2+f_4+...+f_{2m-2})]+O(h^4)\\
    & +\dfrac{h}{2}(f_{2m}+f_{2m+1})+O(h^3)
    +\dfrac{h}{2}(f_{2m+1}+f_{2m+2})+O(h^3)\\
    & +\dfrac{h}{3}[ f_{2m+2}+f_{2n}+4(f_{2m+3}+f_{2m+5}+...+f_{2n-1})+2(f_{2m+4}+f_{2m+6}+...+f_{2n-2})]+O(h^4)\\
    &= \dfrac{h}{3}[ f_0+f_{2n}+4(f_1+f_3+...+f_{2n-1})+2(f_2+f_4+...+f_{2n-2})]\\
    &+ \dfrac{h}{6}(f_{2m}-2f_{2m+1}+f_{2m+2}) +O(h^3)\\
    & = F(G_l (t_i,.)\varphi(.)) +O(h^3)
\end{align*}
Thus, in the both cases of $i$, even or odd, we have the estimate \eqref{eq:Sim1}.\\
The estimate \eqref{eq:Sim2} is obtained analogously as \eqref{eq:Sim1} if taking into account that
$$ 2G_2^* (t_i,t_i)= G_2^- (t_i,t_i)+G_2^+ (t_i,t_i),
$$
where 
\begin{align*}
G_2^{\pm} (t_i,t_i)=\lim _{s\rightarrow t_i \pm 0}G_2(t_i,s)
\end{align*}
\end{Proof}

\begin{Theorem}\label{thm5}
Under the assumptions of Proposition \ref{prop1New}, for the approximate solution of the problem \eqref{eq2} obtained by the discrete iterative method 2 on the uniform grid with gridsize $h$ there hold the estimates
\begin{equation*}
\begin{split}
\|U_k-u\| &\leq \left(M_0+\frac{1}{r}\right)p_kd +O(h^3), \; \|Y_k-u'\| \leq M_1p_kd +O(h^3), \\
\|Z_k-u''\| &\leq M_2p_kd +O(h^3).
\end{split}
\end{equation*}
\end{Theorem}

\section{Examples}
Consider some examples for confirming the validity of the obtained theoretical results and the efficiency of the proposed iterative method. \\

\noindent \textbf{Example 1. } (Problem 2 in \cite{Pandey1}) \\
Consider the problem 
\begin{align*}
\begin{split}
u'''(x)&=x^4u(x)-u^2(x)+f(x), \; 0<x<1,\\
u(0)&=0, \; u'(0)=-1,\; u'(1)=\sin (1),
\end{split}
\end{align*}
where $f(x)$ is calculated so that the exact solution of the problem is $$u^*(x)=(x-1)\sin (x).$$
It is easy to verify that with $M=7$ all conditions of Theorem \ref{theorem3} are satisfied, so the problem has a unique solution.
The results of the numerical experiments with two different tolerances are given in Tables \ref{table:1}- \ref{table:3}.

\begin{table}[h!]
\centering
\caption{The convergence in Example 1 for $TOL= 10^{-4}$}
\label{table:1}
\begin{tabular}{cccccc}
\hline 
$N$ & $K$& $Error_{trap}$ & $Order$ & $Error_{Simp}$ & $Order$\\ 
\hline 
8 & 3 & 9.9153e-04 &  & 9.7143e-04 &  \\ 
16 & 3 & 2.4646e-04 & 2.0083 & 1.3101e-04 & 2.8905  \\ 
32 & 3 & 6.0906e-05 & 2.0167 & 1.6020e-05& 3.0317 \\ 
64 & 3 & 1.4563e-05 & 2.0643 & 1.2587e-06& 3.6696 \\
128 & 3 & 2.9796e-06 & 2.2891 & 8.8553e-07& 0.5073 \\
256 & 3 & 4.3187e-07 & 2.7865 & 8.8165e-07& 0.0063 \\
512 & 3 & 6.7435e-07 & -0.6429 & 8.8118e-07& -7.7719e-04 \\
1024 & 3 & 8.2295e-07 & -0.2873 & 8.8112e-07& 9.6181e-05 \\
\hline 
\end{tabular} 
\end{table}

\begin{table}[h!]
\centering
\caption{The convergence in Example 1 for $TOL= 10^{-6}$}
\label{table:2}
\begin{tabular}{cccccc}
\hline 
$N$ & $K$& $Error_{trap}$ & $Order$ & $Error_{Simp}$ & $Order$\\ 
\hline 
8 & 4 & 9.99237e-04 &  & 9.7223e-04 &  \\ 

16 & 4 & 2.4734e-04 & 2.0044 & 1.3189e-04 & 2.8820  \\ 

32 & 4 & 6.1802e-05 & 2.0008 & 1.6915e-05& 2.9629 \\ 
64 & 4 & 1.5462e-05 & 1.9989 & 2.1492e-06& 2.9765 \\
128 & 4 & 3.8797e-06 & 1.9947 & 2.8688e-07& 2.9053 \\
256 & 4 & 9.8437e-07 & 1.9787 & 5.2749e-08& 2.4439 \\
512 & 4 & 2.6054e-07 & 1.9177 & 2.3446e-08& 1.1698 \\
1024 & 4 & 7.9583e-08 & 1.7110 & 1.9786e-08& 0.2448 \\
\hline 
\end{tabular} 
\end{table}

\begin{table}[h!]
\centering
\caption{The convergence in Example 1 for $TOL= 10^{-10}$}
\label{table:3}
\begin{tabular}{cccccc}
\hline 
$N$ & $K$& $Error_{trap}$ & $Order$ & $Error_{Simp}$ & $Order$\\ 
\hline 
8 & 7 & 9.9235e-04 &  & 9.7222e-04 &  \\ 

16 & 7 & 2.4732e-04 & 2.0045 & 1.3187e-04 & 2.8822  \\ 

32 & 7 & 6.1782e-05 & 2.0011 & 1.6896e-05& 2.9643 \\ 
64 & 7 & 1.5443e-05 & 2.0003 & 2.1301e-06& 2.9877 \\
128 & 7 & 3.8605e-06 & 2.0001 & 2.6774e-07& 2.9923 \\
256 & 7 & 9.6511e-07 & 2.0000 & 3.3544e-08& 2.9965 \\
512 & 7 & 2.4128e-07 & 2.0000 & 4.1977e-09& 2.9984 \\
1024 & 7 & 6.0319e-08 & 2.0000 & 5.2483e-10& 2.9997 \\
\hline 
\end{tabular} 
\end{table}
In the above tables $N$ is the number of grid points, $K$ is the number of iterations, $Error_{trap},\; Error_{Simp} $ are errors $ \| U_K-u^* \|$ in the cases of using Method 1 and Method 2, respectively,$Order$ is the order of convergence calculated by the formula
$$ Order=\log _2 \frac{\|U^{N/2}_K-u^*\|}{\|U^{N}_K-u^*\|}.
$$
In the above formula the superscripts $N/2$ and $N$ of $U_K$ mean that $U_K$ is computed on the grid with the corresponding number of grid points. \\
From the tables we observe that for each tolerance the number of iterations is constant and the errors of the approximate solution decrease with the rate (or order) close to 2 for Method 1 and close to 3 for Method 2 until they cannot improved. This can be explained as follows. Since the total error of the actual approximate solution consists of two terms: the error of the iterative method on continuous level and the error of numerical integration at each iteration, when these errors are balanced, the further increase of number of grid points $N$(or equivalently, the decrease of grid size $h$) cannot in general improve the  accuracy of approximate solution.
 
Notice that in \cite{Pandey1} the author used  Newton-Raphson iteration method to solve  nonlinear system of equations arisen after discretization of the differential problem. Iteration process is continued until the maximum difference between two successive iterations , i.e., $\|U_{k+1}-U_k  \|$ is less than $10^{-10}$. The number of iterations for achieving this tolerance is not reported. The accuracy for some different $N$ is (see \cite[Table 2]{Pandey1})
\begin{table}[h!]
\centering
\caption{The results in \cite{Pandey1} for the problem in Example 1}
\label{table:4}
\begin{tabular}{ccccc}
\hline 
$N$ & 8 & 16 & 32 & 64 \\ 
\hline 
Error & 0.11921225e-01 & 0.33391170e-02 & 0.87742222e-03 & 0.23732412e-03 \\ 
\hline 
\end{tabular} 
\end{table}\\
From the tables of our results and of Pandey it is clear that our method gives much better accuracy.\\

\noindent \textbf{Example 2. } (Problem 2 in \cite{Pandey2}) \\ 
Consider the problem 
\begin{align*}
\begin{split}
u'''(x)&=-xu''(x)-6x^2+3x-6, \; 0<x<1,\\
u(0)&=0, \; u'(0)=0,\; u'(1)=0.
\end{split}
\end{align*}
It is easy to verify that with $M=9$ all conditions of Theorem \ref{theorem3} are satisfied, so the problem has a unique solution. This solution is $u(x)=x^2 (\frac{3}{2}-x)$.
The results of the numerical experiments with different tolerances are given in Tables \ref{table:5}, \ref{table:6} and \ref{table:7}.\\

\begin{table}[h!]
\centering
\caption{The convergence in Example 2 for $TOL= 10^{-4}$}
\label{table:5}
\begin{tabular}{cccccc}
\hline 
$N$ & $K$& $Error_{trap}$ & $Order$ & $Error_{Simp}$ & $Order$\\ 
\hline 
8 & 6 & 0.0078&  & 9.7662e-04 &  \\ 
16 & 6 & 0.0020 & 2.0000 & 1.2215e-04 & 2.9991  \\ 
32 & 6 & 4.8837e-04 & 1.9998 & 1.5345e-05& 2.9929 \\ 
64 & 6 & 1.2216e-04 & 1.9992 & 1.9936e-06& 2.9443 \\
128 & 6 & 3.0604e-05 & 1.9969 & 3.2471e-07& 2.6181 \\
256 & 6& 7.7157e-06 & 1.9878 & 1.1612e-07& 1.4835 \\
512 & 6 & 1.9937e-06 & 1.9524 & 9.0051e-08& 0.3868 \\
1024 & 6 & 5.6316e-07 & 1.8238 & 8.6794e-08& 0.0532 \\
\hline 
\end{tabular} 
\end{table}

\begin{table}[h!]
\centering
\caption{The convergence in Example 2 for $TOL= 10^{-6}$}
\label{table:6}
\begin{tabular}{cccccc}
\hline 
$N$ & $K$& $Error_{trap}$ & $Order$ & $Error_{Simp}$ & $Order$\\ 
\hline 
8 & 8 & 0.0078&  & 9.7662e-04 &  \\ 
16 & 6 & 0.0020 & 2.0000 & 1.2215e-04 & 2.9991  \\ 
32 & 6 & 4.8837e-04 & 1.9998 & 1.5345e-05& 2.9929 \\ 
64 & 6 & 1.2216e-04 & 1.9992 & 1.9936e-06& 2.9443 \\
128 & 6 & 3.0604e-05 & 1.9969 & 3.2471e-07& 2.6181 \\
256 & 6& 7.7157e-06 & 1.9878 & 1.1612e-07& 1.4835 \\
512 & 6 & 1.9937e-06 & 1.9524 & 9.0051e-08& 0.3868 \\
1024 & 6 & 5.6316e-07 & 1.8238 & 8.6794e-08& 0.0532 \\
\hline 
\end{tabular} 
\end{table}

\begin{table}[h!]
\centering
\caption{The convergence in Example 2 for $TOL= 10^{-10}$}
\label{table:7}
\begin{tabular}{cccccccc}
\hline 
$N$ & $K$& $Error_{trap}$ & $Error_{Simp}$ & $N$ & $K$& $Error_{trap}$ & $Error_{Simp}$ \\ 
\hline 
8 & 11 & 0.0078 & 2.0650e-13 & 64 & 11 & 1.2207e-04 & 2.5890e-13\\ 

16 & 11 & 0.0020 & 2.6790e-13 & 128 & 11 & 3.0518e-05 & 2.5790e-13 \\ 

32 & 11 & 4.8828e-04 & 2.6279e-13 & 256 & 11 & 7.6294e-06 & 2.5802e-13\\ 
\hline 
\end{tabular} 
\end{table}

Notice that \cite{Pandey2} the author used  Gauss-Seidel iteration method to solve  linear system of equations arisen after discretization of the differential problem. Iteration process is continued until the maximum difference between two successive iterations , i.e., $\|U_{k+1}-U_k  \|$ is less than $10^{-10}$. The results for some different $N$ are
\begin{table}[h!]
\centering
\caption{The results in \cite{Pandey2} for the problem in Example 2}
\label{table:8}
\begin{tabular}{ccccc}
\hline 
$N$ & 128 & 256 & 512 & 1024 \\ 
\hline 
Error & 0.30696392e-4 & 0.61094761(-5) & 0.14379621e-5 & 0.41723251e-6 \\ 
\hline 
Iter & 53 & 5 & 3 & 4\\
\hline 
\end{tabular} 
\end{table}\\
From the tables of our results and of Pandey it is clear that our method gives better accuracy and requires less computational work.\\

\noindent \textbf{Example 3. } \\ 
Consider the problem 
\begin{align*}
\begin{split}
u'''(x)&=(u(x))^2+u'(x)-e^{2x}, \; 0<x<1,\\
u(0)&=1, \; u'(0)=1,\; u'(1)=e.
\end{split}
\end{align*}
It is easy to verify that with $M=10$ all conditions of Theorem \ref{theorem3} are satisfied, so the problem has a unique solution. This solution is $u(x)=e^x$.
The results of the numerical experiments with different tolerances are given in Tables \ref{table:9} and \ref{table:10}.\\
\begin{table}[h!]
\centering
\caption{The convergence in Example 3 for $TOL= 10^{-4}$}
\label{table:9}
\begin{tabular}{cccccccc}
\hline 
$N$ & $K$& $Error_{trap}$ & $Error_{Simp}$ & $N$ & $K$& $Error_{trap}$ & $Error_{Simp}$ \\ 
\hline 
16 & 8 & 5.4059e-04 & 5.2038e-05 & 128 & 8 & 1.0341e-05 & 2.6902e-06\\ 

32 & 8 & 1.3655e-04 & 1.4204e-05 & 256 & 8 & 4.0312e-06 & 2.1184e-06 \\ 

64 & 8 & 3.5582e-05 & 4.9811e-06 & 512 & 8 & 2.4537e-06 & 1.9755e-06\\ 

\hline 
\end{tabular} 
\end{table}

\begin{table}[h!]
\centering
\caption{The convergence in Example 3 for $TOL= 10^{-6}$}
\label{table:10}
\begin{tabular}{cccccccc}
\hline 
$N$ & $K$& $Error_{trap}$ & $Error_{Simp}$ & $N$ & $K$& $Error_{trap}$ & $Error_{Simp}$ \\ 
\hline 
16 & 11 & 5.3866e-04 & 5.0053e-05 & 128 & 11 & 8.3853e-06 & 7.3348e-07\\ 

32 & 11 & 1.3460e-04 & 1.2241e-05 & 256 & 11 & 2.0750e-06 & 1.6199e-07 \\ 

64 & 11 & 3.3627e-05 & 3.0231e-06 & 512 & 11 & 4.9743e-07 & 1.9180e-08\\ 

\hline 
\end{tabular}
\end{table}

\begin{figure}[ht]
\begin{center}
\includegraphics[height=5cm,width=8cm]{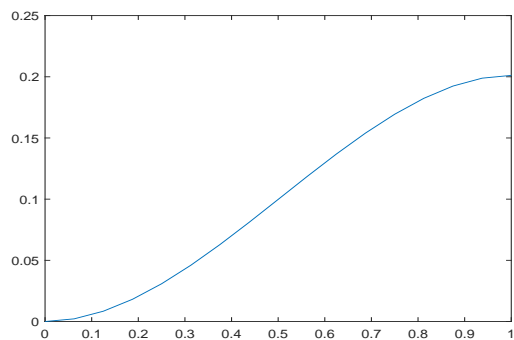}
\caption{The graph of the approximate solution in Example $4$. }
\label{fig3}
\end{center}
\end{figure}

\noindent \textbf{Example 4. }\\ 
Consider the problem for fully third order differential equation
\begin{align*}
\begin{split}
u'''(x)&=-e^{u(x)}-e^{u'(x)}-\frac{1}{10}(u''(x))^2, \; 0<x<1,\\
u(0)&=0, \; u'(0)=0,\; u'(1)=0.
\end{split}
\end{align*}
It is easy to verify that with $M=3$ all conditions of Theorem \ref{theorem3} are satisfied, so the problem has a unique solution.
\begin{table}[h!]
\centering
\caption{The convergence in Example 4 for $TOL= 10^{-10}$}
\begin{tabular}{ccccc}
\hline 
$N$ & 8 & 16 & 32 & 64 \\ 
\hline 
$K$ & 15 & 15 & 15 & 15\\
\hline 
\end{tabular} 
\end{table}\\
The numerical solution of the problem is depicted in Figure \ref{fig3}.

\section{Conclusion}
In this paper we  established the existence and uniqueness of solution for a boundary value problem for fully third order differential equation. Next, for finding this solution we proposed an iterative method at both continuous and discrete levels. The numerical realization of the discrete iterative method is very simple. It is based on popular rules for numerical integration.
One of the important results is that we obtained an estimate for total error of the approximate solution which is actually obtained. This total error depends on the number of iterations performed and the discretization parameter. The validity of the theoretical results and the efficiency of the iterative method are illustrated in examples. \par
The method for investigating the existence and uniqueness of solution and the iterative schemes for finding solution in this paper can be applied to other third order nonlinear  boundary value problems, and in general, for higher order nonlinear  boundary value problems.


\end{document}